\documentclass[12pt]{amsart}
\usepackage{epsfig}
\usepackage{amsmath}
\usepackage{amsfonts}
\usepackage{amssymb}
\usepackage{epigraph, fancyhdr}

\def\K{\mathcal K}
\def\J{\mathcal J}

\def\O{\mathcal O}

\def\1{\mathbf 1}
\def\M{{\overline{\mathcal M}}}
\def\N{{\mathcal N}}

\def\QQ{\mathbb Q}
\def\ZZ{\mathbb Z}
\def\CC{\mathbb C}

\def\Res{\operatorname{Res}}

\def\det{\operatorname{det}}

\def\tilde{\widetilde}

\def\t{{\mathbf t}}

\def\gs{\sigma}
\def\gl{\nu}
\def\gL{\Lambda}

\def\lan{\langle}
\def\ran{\rangle}

\def\tr{\operatorname{tr}}

\renewcommand{\Delta}{\triangle}

\title[Fixed point localization]
      {Permutation-equivariant \\ quantum K-theory II. \\
      Fixed point localizaion}

\author[A. Givental]{Alexander GIVENTAL}
\thanks{This material is based upon work supported by the National 
Science Foundation under Grant DMS-1007164, and by the IBS Center for Geometry 
and Physics, POSTECH, Korea.}

\begin{document}

\begin{abstract}
  
On projective spaces as examples of toric manifolds, we examine K-theoretic fixed point localization. On the one hand, we will see how the permutation-equivariant theory of the point target space emerges as a necessary ingredient. On the other hand, we will completely characterize
the genus-0 permutation-equivariant quantum K-theory of the given toric manifold in terms of such theory for the point, and a certain recursion relation.  

\end{abstract}

\date{June 31, 2015}

\maketitle

\section*{Example: Localization on $\M_{0,1}(\CC P^1, 2)$}

Let the target space $X$ be $\CC P^N=\operatorname{proj}(\CC^{N+1})$. We consider it as a toric manifold equipped with the action of the torus $T^{N+1}$
of matrices $\operatorname{diag}(\gL_0,\dots,\gL_N)$ acting naturally on $\CC^{N+1}$. One can define K-theoretic GW-invariants, ordinary or permutation-equivariant, {\em equivariant with respect to the torus action} and taking values in the ring $Repr (T^{N+1})=\ZZ [\gL_0^{\pm 1},\dots,\gl_N^{\pm 1}] $. It is imperative in this case to extend our ground ring $\gL$ by including Novikov's variable(s), $Q$, as well as $\gL_i^{\pm 1}$, and in the permutation-equivariant case, extend the structure of $\lambda$-algebra by $\Psi^m(Q^d)=Q^{md}$ and $\Psi^k(\gL_i)=\gL_i^m$.

In fact, for complex projective spaces, one value of the small J-function has been known for a long time even in the $T$-equivariant setting:
\begin{align*} \J_{\CC P^N}(0)&:=1-q+\sum_{i,d}\phi_iQ^d\lan \frac{\phi^i}{1-qL}\ran_{0,1,d} =\\
 &(1-q)\sum_{d\geq 0} \frac{Q^d}{\prod_{i=0}^N(1-qP\gL_i^{-1})(1-q^2P\gL_i^{-1}) \cdots (1-q^dP\gL_i^{-1})}.\end{align*}
  It takes values in the torus-equivariant K-ring of $\CC P^N$, which
  is generated over $\gL$ by the Hopf bundle $P$ satisfying the relation
  \[ (1-P\gL_0^{-1})\cdots (1-P\gL_N^{-1})=0. \]
Since the input of the J-function is $0$, it belong to the permutation-equivariant theory and the non-equivariant theory as well. In the latter capacity it was computed {\em ad hoc} in \cite{GiL} (on the basis of ``linear sigma-models'' and some rationality arguments)

The specializations to one of the roots $P=\gL_i$ represent the restrictions to the $N+1$ fixed points of the torus on $\CC P^N$, i.e. the components in
the basis  $\{ \phi_i\} $ of ``delta-function'' of the fixed points.
Say, for $i=0$, we have
  \[ J_{\CC P^N}(0)_{P=\gL_0} = (1-q)\sum_{d\geq 0}\frac{Q^d}{\prod_{r=1}^d(1-q^r) \prod_{i=1}^N\prod_{r=1}^d(1-q^r\gL_0\gL_i^{-1})}.\]
Our aim is to develop  a systematic way to characterize K-theoretic GW-invariants, albeit in genus-$0$ for the moment, using localization to fixed point of $T$-action on moduli spaces of stable maps. But first we would like to see how the above expression can emerge in the simplest example of $N=1$ and degree $d\leq 2$.  

For $d=0$, we have $1-q$, the ``dilaton shift'' summand.

For $d=1$, we have $1/(1-q\gL_0\gL_1^{-1})$, which makes sense, because this is the value of $1/(1-qL)$ at this fixed point. Indeed, the moduli space of degree-$1$ rational
curves in $\CC P^1$ with $1$ marked point is $\CC P^1$ itself: there is only one
such map, the identity map, and $\CC P^1$ worth of choices for the marked
point. When the marked point is localized to the $T$-fixed point (by the choice
of $\phi^i$ with $i=0$), the line ``bundle'' $L$ becomes the cotangent line to
$\CC P^1$ at this fixed point, and $\gL_0\gL_1^{-1}$ is the character by which the torus acts on this line.

For $d=2$, we have a reduced rational function in $q$ with denominator of degree $5$, which we, abbreviating $\gL_0\gL_1^{-1}$ to $\lambda$,  decompose into elementary fractions: 
\begin{align*} \frac{1}{(1-q^2)(1-\lambda q)(1-\lambda q^2)}=&
 \frac{1}{2(1-\lambda)^2(1-q)}+\frac{1}{2(1-\lambda^2)(1+q)}\\
+\frac{\lambda^3}{(1-\lambda)(1-\lambda^2)(1-\lambda q)}&-
\sum_{\pm}\frac{\lambda}{2(1-\lambda)(1\pm \sqrt{\lambda})(1\pm \sqrt{\lambda}q)}.\end{align*} 
Note that the rational function, expanded in a $q$-series, will have coefficients which are Laurent polynomials in $\lambda$ -- characters of representation so the torus on cohomology of the sheaf $\phi^0 L^m$ on the moduli space $\M_{0,1}(\CC P^1, 2)$. However, the elementary fractions have poles at $\lambda=\pm 1$. This suggest that they make sense in the context
fixed point localization. Yet, we also need to understand the occurrence of
$\sqrt{\lambda}$ which makes sense only on the double-cover of $T^2$.

Recall that Lefschetz' fixed point formula on a com[act complex manifold
 $M$ equipped with a holomorphic action of a torus $T$ and an equivariant holomorphic bundle $V$ computes the supertrace 
  $\tr_{\lambda} H^*(M; V)$ for $\lambda\in T$, i.e. the character of $T$ on the sheaf cohomology of $V$, as a holomorphic Euler characteristic $\chi$ on the fixed point manifold $M^T$:
  \[ \tr_h H^*(M ; V) = \chi \left(M^T; \frac{\tr_{\lambda} V}{\tr_{\lambda}
    \bigwedge^* \N^*_{M^T}}\right).\]
Here $\N^*_{M^T}$ is the conormal bundle to the fixed point manifold, and the trace of a vector bundle restricted to the fixed point locus is a virtual bundle defined by decomposing the bundle according to eigenvalues  $\lambda_j$:
    \[ \tr_{\lambda} V := \oplus_j \lambda_j V(\lambda_j).\]
In the case when all fixed points are isolated we obtain the sum
    \[ \tr_{\lambda} H^*(M; V) = \sum_{p\in M^T} \frac{\tr_{\lambda} V_p} {\det_{1-\lambda}T^*_pM }.\]
    The denominators here have poles at some points on of finite order on $T$ --- but no square roots. For example, the fixed point formula on $\CC P^1$
    looks this way:
    \begin{align*} \tr_{(\gL_0,\gL_1)}&H^*(\CC P^1; P^k) =\frac{\gL_0^k}{1-\gL_0\gL_1^{-1}} + \frac{\gL_1^k}{1-\gL_1\gL_0^{-1}}=\\
 &-\text{\large$\Res_{P=0,\infty}$}\left( \frac{P^k}{(1-P\gL_0^{-1})(1-P\gL_1^{-1})}\frac{dP}{P}\right) .
\end{align*}    
    
    To understand the origin of $\sqrt{\lambda}$ in the earlier formula, we need to realize that the moduli $\M_{0,1}(\CC P^1, 2)$ is an orbifold. Of course, the torus $T^2$ acting on $\CC P^1$, also transforms stable maps to $\CC P^1$, and thereby acts on the
    space of their equivalence classes. However the cotangent spaces $T*_p\M$
    at the fixed points (and more generally the conormal bundle $\N^*_{\M^T}$ to the fixed point suborbifold) at the cotangent spaces (conormals) on the orbifold chart, while the orbifold locally is the quotient of this chart
    by the local symmetry group $\Gamma_p$. The infinitesimal action of the torus on $\M$ lifts to the chart unambiguously, but integrating it to the group action results in its extension by $\Gamma_p$. The extension may be
    a connected covering of the torus (as the double-cover $\pm \sqrt{\lambda}$ above), or the product $T\times \Gamma_p$ (or, of course, any combination of these possibilities). In any case to compute $\tr_{\lambda} H^*(\M ; V)$, one needs to take in Lefschetz' formula the {\em average} over all $|\Gamma_p|$ preimages of $\lambda \in T$ in the extension group. Having all this in mind, let us explicitly apply Lefschetz' formula on $\M_{0,1}(\CC P^2, 2)$.

    The following diagram accounts for all fixed points of of the torus action
    on th moduli space of degree-$2$ stable maps to $\CC P^1$ with one marked point, mapped to one of the two fixed point on $\CC P^1$, namely the one with the index $i=0$.

\begin{figure}[htb]
\begin{center}
\epsfig{file=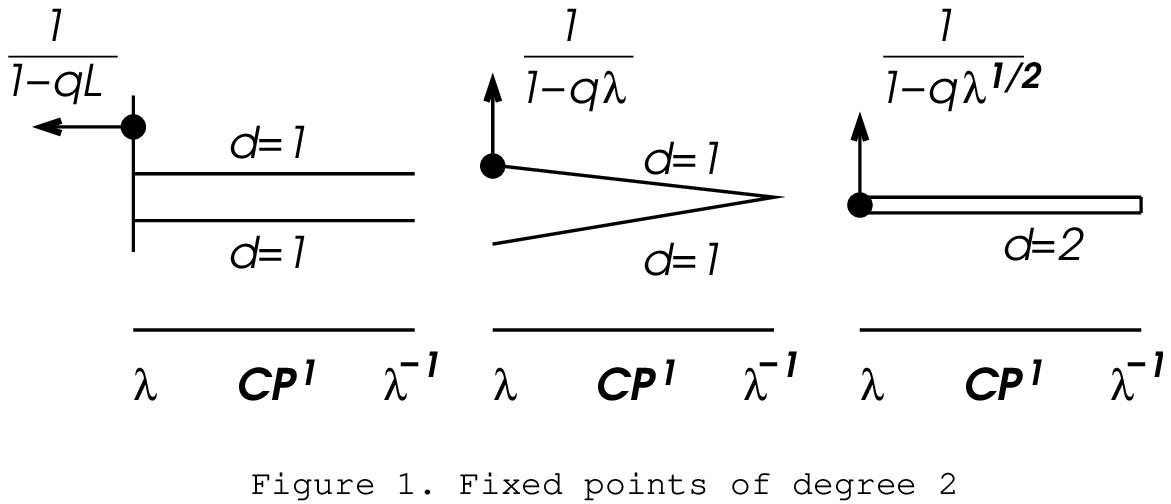}
\end{center}
\end{figure}    
    
The bottom intervals represent the target $\CC P^1$ shown together with the eigenvalues $\lambda^{\pm1}$ on the cotangent lines at the fixed points (the ends of the interval). The top pictures represent possible nodal rational curves projected to $\CC P^1$ in such a way that their equivalence classes
are torus-invariant. The thick dot represents the marked point and must be mapped to the ``left'' fixed point in $\CC P^1$ to contribute non-trivially
to the correlator $\lan \phi^0/(1-qL)\ran_{0,1,2}$, for $\phi^0=0$ at the ``right'' fixed point. At the ``left'' fixed point $\phi^0=1-\lambda$, as it is evident from the above localization formula on $\CC P^1$ and $\phi_0=1$.

The leftmost picture represents a $\CC P^1$ mapped to the ``left'' fixed point,
carrying the marked point, and two nodes where two copies of $\CC P^1$, mapped with degree $d=1$ to the target $\CC P^1$, are attached.

The middle picture represents a $\CC P^1$ carrying the marked point
and mapped with degree $d=1$ to the target, with another degree-$1$ copy
of $\CC P^1$ attached to it at the node over the ``right'' fixed point.

The rightmost picture represents a degree $d=2$ map $\CC P^1\to \CC P^1: z\mapsto z^2$, ramified at $z=0$ and $\infty$ over the ``left'' and ``right''
fixed points, and carrying the marked point at one of them. The arrow symbolizes the {\em horn} of the curve, by which we call the marked point carrying the input $1/(1-qL)$. This notion, currently redundant, will become
useful afterwards.

We can now match contributions of the fixed points with the above elementary fractions.

The leftmost curve has a $\ZZ_2$-symmetry interchanging the degree-$1$ branches. It acts on the moduli space $\M_{0,3}$ of ``vertical'' components
  by interchanging the nodes. Of course, $\M_{0,3}$ is a point, but the action
  is non-trivial on the cotangent line $L$ at the marked point, and on the
  conormal space to the fixed point. The two elements of $\ZZ_2$ contribute
  the elementary fraction $1/(1\mp q)$ (corresponding to the eigenvalues $L=\pm 1$ respectively on the cotangent line). The moduli space has dimension
  $3$. The cotangent space consists of three $1$-dimensional modes of deformation: the shift of the vertical component away from the fixed point in the target $\CC P^1$, and two modes of ``smoothing the nodes'' (i.e. turning the nodal $xy=0$ into $xy=\epsilon$). The smoothing of the two nodes can be performed $\ZZ_2$-symmetrically or anti-symmetrically. This shows that the determinants on the cotangent $3$-space, occurring in the denominator of Lefschetz' formula, are equal to $(1-\lambda)^3$ for the identity element of $\ZZ_2$, and $(1-\lambda)^2(1+\lambda)$ for the non-identity one. Yet one of the factors $(1-\lambda)$ cancels with $\phi_0$. The factor $1/2$ stands for taking the average over $\ZZ_2$. The outcomes match exactly the first two elementary fractions. This analysis shows that {\em fixed point localization in quantum K-theory of a target space $X$ engages  permutation-equivariant quantum K-theory of the fixed point manifold $X^T$, i.e. if the point target space if the $X^T$ is discrete.}

  The middle curve contributes the elementary fraction $1/(1-\lambda q)$ (with the eigenvalue on the cotangent line $L$ at the marked point equal $\lambda$) while the rightmost double-cover curve brings in the average of two fractions $(1\pm \sqrt{\lambda}$. We leave it for the reader to verify that the coefficients at the elementary fractions agree with the determinants on the conormal $3$ space to the corresponding fixed points in the moduli space $\M_{0,1}(\CC P^1, 2)$. 

\section*{Localization as recursion}

Now we systematically explore fixed point localization in genus-$0$ quantum K-theory a target space $X$, equipped with an action of complex torus $T$, assuming that the fixed points are isolated, and one dimensional orbits are
isolated as well (though the latter condition can be relaxed), focusing on the example $X=\CC P^N$ for the sake of simplicity.

The method goes back to th paper \cite{Ko} by M. Kontsevich. The equivalence
of stable maps to $X$ is defined via identifications of the domain curves together with markings. A stable map $(\Sigma, \gs_1m\dots, \gs_n) \to X$, in order to represent a torus-invariant equivalence class, must have a torus-invariant range. Thus, $\Sigma$ must fall apart into irreducible components mapped to $X^T$ (let us call them {\em vertices}) connected by chains of copies of $\CC P^1$ ({\em legs}, see the middle picture on Figure 1 showing such a chain of two legs), each mapped onto the closure of a one one-dimensional orbit. Such an orbit connects two distinct fixed points, and together with them, also forms a $\CC P^1$. Each leg $\CC P^1\to \CC P^1$ is a map $z\mapsto z^m$ of certain {\em multiplicity} $m=1,2,\dots$ ramified at $z=0, \infty$ over the fixed points connected by the orbit. All marked points must either be located on the vertices, or at the end of a chain of legs. 

Now we express the big J-function $\J_X$ of permutation-equivariant quantum K-theory on $X=\CC P^N$ in terms of fixed point localization.

In the basis $\{ \phi \}$ of delta-functions of the fixed points,  we have: $\J_X(\t)= \sum \J_X^{(i)}(\t) \phi_i$, $\t=\sum \t^{(i)} \phi_i$, $i=0,\dots, N$. The diagrams on Figure 2 represent various fixed point contributions to
\[ J_X^{(i)}(\t) = (1-q)+\t^{(i)}(q,q^{-1})+\sum_{n,d}Q^d\lan \t,\dots,\t,
\frac{\phi^i}{1-qL}\ran_{0,n+1,d}.\]

\begin{figure}[htb]
\begin{center}
\epsfig{file=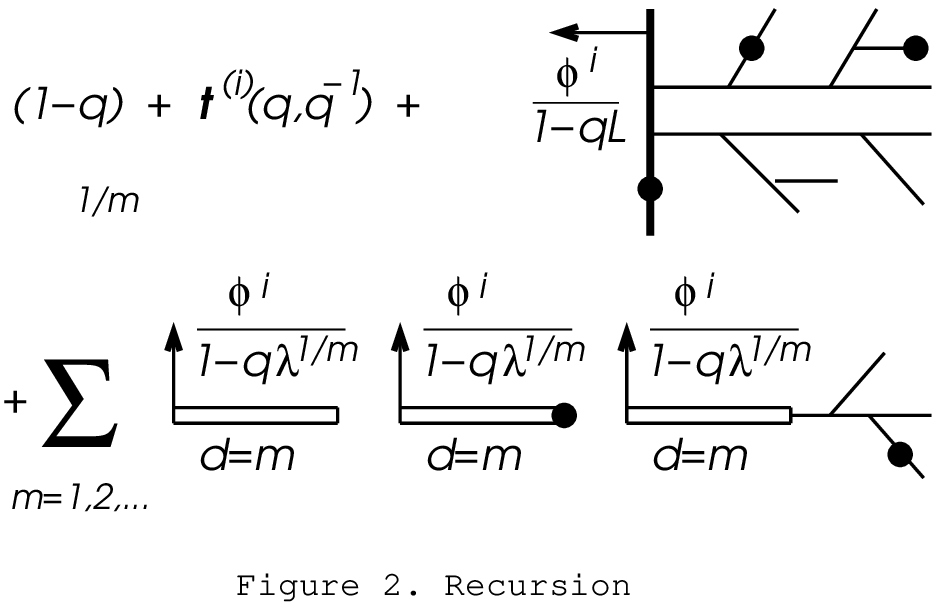}
\end{center}
\end{figure}   

Namely, the marked point carrying the input $\phi^i/(1-qL)$ (the {\em horn}) lies either on a stable component ({\em vertex}) mapped to the $i$-th fixed point in $\CC P^N$ or on a leg of certain multiplicity $m$, issued from the $i$-th fixed point and heading toward a $j$-th fixed point with $j\neq i$
along the straight line in $\CC P^N$ connecting these fixed points.

In the latter case (illustrated by the bottom pictures), the cotangent line $L$ at the horn carries the action of the torus given by the character $\lambda^{1/m}$, where $\lambda=\gL_i \gL_j^{-1}$
(compare with $\gL_0 \gL_1^{-1}$ in the case of $\CC P^1$). As a function of
$q$, such contributions into $\J_X^{(i)}$ form the elementary fraction proportional to $1/(1-\lambda^{1/m} q)$ (to be averaged over all the $m$ values of the $m$-th root of $\lambda$). We postpone the discussion of the proportionality coefficient (which depends on the entire tree attached at the
other end of the leg). 

In the former case (illustrated by the top picture), vertices form moduli
space $\M_{0,r+1}$, where $r$ ($=3$ on the diagram) is equal to the total number of marked points and legs attached to the vertex. As functions of $q$, such  fixed points contributions have poles at roots of unity. Indeed, as in the example with degree $2$ curves in $\CC P^1$, in Lefschetz' contributions of the points in $\M_{0,r+1}$ fixed by a discrete symmetry $h\in S_r$, the trace $\tr_h (1/(1-qL)=1/(1-q\zeta \tilde{L})$. Here $\zeta$ is the eigenvalue of $h$
on $L$, and $\tilde{L}$ is the restriction of $L$ to the fixed point locus of $h$. Hence $\zeta$ is a root of unity, and $(\tilde{L}-1)^s=0$ in the K-ring of the fixed point locus of dimension $<s$. Therefore only elementary fractions $1/(1-\zeta q)^k$ with $k\leq s$ can arise.

Note that 
\[ \frac{1}{1-\tilde{L}\mu} = \frac{1}{(1-\mu)-\mu(\tilde{L}-1)}=
\sum_{k\geq 0}\frac{\mu^k(\tilde{L}-1)^k}{(1-\mu)^{k+1}}.\]
Using this with $\mu =\zeta \lambda^{1/m}=\zeta \gL_i^{1/m}\gL_j^{-1/m}$, we find
that fractions of the form $1/(1-q\lambda^{1/m})$ (on the bottom of the diagram) can be considered as legitimate {\em inputs} in
permutation-equivariant GW-invariants of the point target space,
{\em if one localizes the coefficient ring $\gL$ to allow the division by
  $1-\zeta \gL_i^{1/m}\gL_j^{-1/m}$. With such division enabled, we claim that
\begin{align*} \J_X^{(i)}(\t)=\J_{pt}(t)&:=(1-q)+t(q)+
  \sum_{r>1}\lan t(L),\dots,t(L),\frac{1}{1-qL}\ran_{0,r+1}^{S_r},\\
\text{where}\ \ t(q)&:=\t^{(i)}(q,q^{-1})+ \text{bottom part of}\ \J_X^{(i)}(\t).
\end{align*}}
Indeed, contributions of all configurations on the top of the diagram sum up
to the correlator sum of $\J_{pt}(t)$ with the input $t$ which combines all
the remaining terms of the diagram but $1-q$, only with $q$ replaced by $L$. Namely, input $t^{(i)}(L,L^{-1})$ corresponds to a marked point located directly on vertex, while the bottom configurations (with $q$ replaced by $L$) represent all possible branches (on the top) attached to the vertex. In Lefschetz' formula, the factor $1/(1-L\lambda^{1/m})$ represents the conormal
direction to the fixed point locus, which corresponds to smoothing the curve
at the node where the branch meets the vertex. Here the torus acts by the character $\lambda^{1/m}$ on the cotangent line to the leg at the node, while $L$ is the cotangent line bundle to the vertex at this node.

Our next observation characterizes the coefficients at the elementary fractions $1/(1-\lambda^{1/m}q)$. Namely {\em functions $\J_X^{(i)}$, $i=0,\dots, N$, have simple poles at $q=\lambda^{-1/m}$, where $\lambda=\gL_i/\gL_j$, $j\neq i$,  $m=1,2,\dots$, and the roots $\lambda^{1/m}$ assuming all $m$ possible values; the residues at these poles satisfy the following system of relations:
  \[ \Res_{q=\lambda^{-1/m}} \J_X^{(i)}(\t; q) \frac{dq}{q} =-\frac{Q^m}{m}\frac{\phi^i}{C_{ij}(m)}\, \J_X^{(j)}(\t; \lambda^{-1/m}),\]
where the coefficients $C_{ij}(m)$ are certain rational functions on the torus $T$, which can be explicitly written in terms of the $1$-dimensional orbit connecting the $i$-th and $j$-th fixed point and multiplicity $m$.}  

Indeed, the bottom picture on the diagram represent possibilities one encounters with fixed point contributions having the pole at $q=\lambda^{-1/m}$.

It could be a leg of multiplicity $m$ (connecting the $i$-th and $j$-th fixed points), which has a marked point at the other end (the middle picture). This is an isolated fixed point $p$ in the moduli space $\M=\M_{0,2}(\CC P^N, m)$. 
The contribution will have the form
\[ \frac{\phi^i}{(1-\lambda^{1/m}q)}\,\frac{Q^m}{m C_{ij}(m)}\,\t^{(j)}(\lambda^{-1/m},\lambda^{1/m}),\]
where $C_{ij}(m) =\tr_T \bigwedge^* T^*p\M$ is Lefschetz' denominator at the fixed point in question, $\phi_i=\prod_{j\neq i}(1-\gL_i\gL_j^{-1})$ is
such a denominator at the $i$-th fixed point on the target space, and the extra factor $1/m$ weighs the individual contributions in the average over the $m$ values of $\lambda^{1/m}$.

It could be such a leg with a regular (unmarked) point at the other end
(shown on the bottom left picture). This is an isolated fixed point in $\M_{0,1}(\CC P^N, m)$, with the cotangent space ``smaller'' than in the previous case by the $1$-dimensional summand, the cotangent line to the leg at the right end. Consequently the contribution
reads:
\[ \frac{\phi^i}{(1-\lambda^{1/m}q)}\,\frac{Q^m}{C_{ij}(m)}\, (1-\lambda^{-1/m}),\]

Finally, it could be a leg with a {\em tail} attached at the other end (as illustrated on the bottom right picture). The contribution factors as
\[ \frac{\phi^i}{(1-\lambda^{1/m}q)}\,\frac{Q^m}{C_{ij}(m)} \times \text{(tail contribution)}.\]
The smoothing factor at the node, where the leg and tail meet, has the form $1/(1-\lambda^{-1/m}L)$, where $L$ is the cotangent line to the tail at the node. This factor plays the role of the horn for the tail curve, which can be any $T$-fixed point with the horn mapped to the $j$-th fixed point. The total contribution of such tails is the fixed point expression for
\[ \J_X^{(j)}(\t, q)-(1-q)-\t^{(j)}(q,q^{-1}) \]
with $q$ replaced by $\lambda^{-1/m}$.

Collecting all the three cases together yields what we promised.

\medskip

We have shown that for any input $\t$, the components $f^{(i)}:=\J_X^{(i)}(\t)$ of the big J-function in the permutation- and torus-equivariant quantum K-theory of $X=\CC P^N$, when considered as a (series of in Novikov's variables) of rational functions in $q$, satisfy two local requirements:

\smallskip

{\bf (i)} when expanded as meromorphic functions with poles only at the roots of unity, $f^{(i)}$ represent values of the big J-function $\J_{pt}$ of permutation-equivariant quantum K-theory of the point target space;

{\bf (ii)} outside $q=0,\infty$ and the roots of unity, $f^{(i)}$ may have first order poles only at $q=(\gL_j/\gL_i)^{1/m}$, $m=1,2,\dots$, and their residues at these poles satisfy 
\[ \Res_{q=(\gL_j\gL_i)^{1/m}}  f^{(i)}(q) \frac{dq}{q} = -\frac{Q^m}{C_{ij}(m)}
f^{(j)}((\gL_j/\gL_i)^{1/m}),\]
where $C_{ij}(m)$ are rational functions of $(\gL_0,\dots, \gL_N)$ (implicitly
described above explicitly below).

\smallskip

We claim that, conversely, if $f^{(i)}$ satisfy the requirements (i) and (ii), then $\sum_i f_i^{(i)}(q)$ represents a value of $\J_X$.

To be more precise, let us assume (for the sake of certainty) that the lambda-algebra of symmetric
functions $\QQ[[N_1,N_2,\dots ]]$, where $N_k$ are the Newton polynomials
(in infinitely or finitely many variables), is chosen to define
permutation-equivariant correlators. It has been extended by Novikov's
variable(s) $Q$ and by Laurent polynomials in $(\gL_0,\dots,\gL_N)$ to deal with $T$-equivariant quantum K-theory:
\[ \gL = \QQ [[N_1,\N_2,\dots ]] [[Q]] [\gL_0^{\pm 1},\dots, \gL_N^{\pm 1}].\]
We assume that $f^{(i)}$ here are series in $N_k$ and $Q$ with coefficients which are rational functions of $q$ (yet, with coefficients which are Laurent  
polynomials in $\gL_i$. In other words, modulo any power of the ideal generated
by $Q$ and $N_k$, they are rational as functions of $q$. We will also assume that $f^{(i)}=1-q \mod (N_1,N_2,\dots)$, and denote by $\K$ the space of all such functions (which we will nick-name ``rational''). Denote by $\K_{+}$ the subspace of rational functions which are Laurent polynomials in $q$
(modulo any power of that ideal), and by $\K_{+}$ the complementary subspace of
those rational functions which have no pole at $q=0$ and whose difference with
$1-q$ tends to $0$ at $q=\infty$. Given such $f^{(i)}$, $i=0,\dots, N$, we take their Laurent polynomial parts (defined as the projection to $\K_{+}$ along $\K_{-}$) for $(1-q)+\t^{(i)}$, and claim that

\smallskip

{\em the requirements (i) and (ii) together determine $f=\sum f^{(i)}\phi_i$ in terms of $\t:=\sum \t^{(i)}\phi_i$ recursively by degrees of $Q$, and thus identify $f$ with $\J_X(\t)$.}   

\smallskip

Indeed, modulo $Q$, relations (ii) show that $f$ can have poles only at $q=0,\infty$, or the roots of unity. Then (i) does not require any localization, and therefore means that $f^{(i)} \mod Q = \J_{pt}(\t^{(i)}\mod Q)$, which is exactly what $J_X^{(i)}(\t) \mod Q$ are. Now, to determine $f \mod Q^d$, it suffices to reconstruct the part of it in $\K_{-}$, given the Laurent polynomial part $\t \mod Q^d$.

This is done in two steps. First, the elementary fraction in $f^{(i)} \mod Q^d$ with poles away from $q=0,\infty$ and the roots of unity, are recovered from relation (ii), whenever all $f^{(j)}$ with $j\neq i$ are given in degrees $<d$ in $Q$. Then, the expansions of these fractions a power series in $q$ together with $\t^{(i)}\mod Q^d$ form the input in $\J_{pt}$ (see Figure 2), whose value,
according to relation (i), uniquely determines modulo $Q^d$ the part of $f^{(i)}$ with poles at he roots of unity. Thus, we have proved the following

\medskip

{\tt Theorem} (fixed point localization). {\em The range of $\J_X=\sum \J_X^{(i)}\phi_i$ is completely characterized as the set of functions $f=\sum f^{(i)}(q)\phi_i$ such that the rational functions  $f^{(i)}\in \K$ obey local requirements (i) and (ii).}

\medskip

{\tt Remark.} The use of localization for recursive description of $\J_X$
goes back to our paper \cite{GiG}, but the idea of characterizing the vertex
contributions in terms of $\J_{pt}$ was introduced, albeit in the cohomological context, more recently by J. Brown \cite{Br}. After \cite{GiT}, it became clear that this recursion is an instance {\em adelic characterization}, which
will be discussed in Part III.


\section*{Application to $q$-hypergeometric functions}

Trying to apply the localization theorem to certain $q$-hypergeometric series associated with toric manifolds, one immediately finds that the series satisfy the recursive requirement (ii). Let us do this for $X=\CC P^N$. Take
the $q$-hypergeometric series (which as we know from the literature is
$\J_{\CC P^N}(0)$):
\[ J=(1-q)\sum_{d\geq 0}\frac{Q^d}{\prod_{i=0}^N\prod_{r=1}^d(1-P\gL_i^{-1}q^r)}\ \ \ \ \mod \prod_{i=0}^N (1-P\gL_i^{-1}).\]
Then $J^{(i)} = J|_{P=\gL_i}$, $i=0,\dots N$:
\[ J^{(i)}=(1-q)\sum_{d\geq 0}\frac{Q^d}{\left(\prod_{r=1}^d (1-q^r)\right)
  \prod_{j\neq i} \prod_{r=1}^d (1-q^r\gL_i\gL_j^{-1})}.\]
Let us extract from $J^{(0}$ the elementary fraction with the pole at
$q=(\gL_j/\gL_0)^{1/m}$ with $j\neq 0$.

First, from the term with $d=m$, we extract
\[ \frac{(1-(\gL_j/\gL_0)^{1/m})}{(1-q\gL_0/\gL_j)}\, \frac{Q^m}{m}\,\frac{\phi^0}{C_{0,j}(m)},\]
where $\phi_0=\prod_{i\neq 0} (1-\gL_0/\gL_i)$, and
\[ C_{0,j}(m)=\phi^0 \phi^j\,\prod_{r=1}^{m-1}\prod_{i=0}^N\left(1-(\gL_j/\gL_i)^{r/m}
(\gL_0/\gL_i)^{(m-r)/m}\right).\]
The factor $m$ comes here as the limit
\[ \lim_{q\to \lambda^{-1/m}} \frac{1-q^m\lambda}{1-q\lambda^{1/m}},\ \ \ \lambda=\gL_0/\gL_j.\]

The terms of the series $\J^{(0)}_X$ with $d>m$ are products of the term with
$d=m$ and the factor
\[ \frac{Q^{d-m}}{\prod_{r=1}^{d-m}\prod_{i=0}^N (1-q^mq^r\gL_0/\gL_i)}.\]
Replacing here each $q^m$ by $\gL_j/\gL_i$ we obtain the term of $Q$-degree $d-m$ in $\J^{(j)}(q)/(1-q)$. Thus, summing all the elementary fractions
with the pole at $q=(\gL_j/\gL_0)^{1/m}$, we obtain
\[ \frac{1}{1-q(\gL_0/\gL_j)^{1/m}}\,\frac{Q^m}{m}\, \frac{\phi^0}{C_{0,j}(m)}\,
\J_X^{(j)}((\gL_j/\gL_0)^{1/m})\]
as expected.

The product $C_{0,j}(m)$ in the denominator can be identified with the trace $\tr_T$ on the cotangent space to $\M_{0,2}(\CC P^N, m)$ at the appropriate
fixed point. Namely,
as $\CC P^N=\operatorname{proj}(\CC^{N+1})$,
the tangent bundle to $\CC P^N$ can be described $T$-equivariantly as $\CC^{N+1}\otimes \O(1)-\O$. On the $m$-multiple cover of $\CC P^1$ which
is the projectivization of the plane $\CC^2$ with the torus action given
by $(\gL_0, \gL_j)$, the pull-back of $\CC^{n+1}\otimes \O(1)$ has  $(N+1)(m+1)$-dimensional space of sections, which splits into $T$-invariant lines with obvious eigenvalues 
\[ \gL_i \otimes \gL_0^{-r/m}\gL_j^{-s/m}, r,s\geq 0, r+s=m, i=0,\dots, N.\]
The two trivial eigenvalues (they correspond to $(i,r,s)=(0,m,0)$ and $(j,0,m)$) are to be dropped from the list, because they cancel with
the infinitesimal automorphisms of the domain $\CC P^1$ preserving two marked points, with the line of sections of $-\CC$. The rest
gives the list of eigenvalues $\mu$ on the {\em tangent} space to $\M_{0,2}(\CC P^N, m)$ at the point of our interest. Multiplying $(1-\mu^{-1})$ we obtain the
denominator in Lefschetz' formula. 

\medskip

Thus, to qualify for a value of $\J_X$, the $q$-hypergeometric series needs
to pass the test (i), i.e. represent a value of the big J-function $\J_{pt}$
of the permutation-equivariant quantum K-theory of the point, when interpreted as a meromorphic function with poles only at the roots of unity.
In Part III, we will completely describe the range of $\J_{pt}$, and then in Part IV return to the problem of verifying requirement (i).

\enddocument